\documentclass{agtart_a}
\pdfoutput=1

\usepackage{pinlabel,picins}
\usepackage[all]{xy}

\title{Slicing Bing doubles}

\author{David Cimasoni}
\givenname{David}
\surname{Cimasoni}
\address{Department of Mathematics\\UC Berkeley\\970 Evans
Hall\\\newline
Berkeley, CA 94720\\USA}
\email{cimasoni@math.berkeley.edu}
\urladdr{}

\volumenumber{6}
\issuenumber{}
\publicationyear{2006}
\papernumber{83}
\startpage{2395}
\endpage{2415}

\doi{}
\MR{}
\Zbl{}


\received{16 September 2006}
\revised{}
\accepted{18 November 2006}
\published{13 December 2006}
\publishedonline{13 December 2006}
\proposed{}
\seconded{}
\corresponding{}
\editor{CPR}
\version{}

\arxivreference{math.GT/0609458}

\let\xysavmatrix\xymatrix
\def\xymatrix{\disablesubscriptcorrection\xysavmatrix}
\AtBeginDocument{\let\bar\wbar\let\tilde\wtilde\let\hat\what\let\widehat\wwhat\renewcommand{\C}{\mathcal C}}
\makeop{int}

\makeatletter
\def\cnewtheorem#1[#2]#3{\newtheorem{#1}{#3}[section]
\expandafter\let\csname c@#1\endcsname\c@thm}
\makeatother

\newtheorem{thm}{Theorem}[section]
\newtheorem*{thmintro}{Theorem}
\cnewtheorem{lemma}[thm]{Lemma}
\cnewtheorem{prop}[thm]{Proposition}
\cnewtheorem{cor}[thm]{Corollary}
\theoremstyle{remark}
\newtheorem*{rem}{Remark}
\newtheorem*{qu}{Question}
 
\newenvironment{romanlist}
        {\begin{enumerate}
        \renewcommand{\theenumi}{\roman{enumi}}}
        {\end{enumerate}}

\newcommand{\Cc}{\mathbb C}
\newcommand{\N}{\mathcal N}
\newcommand{\A}{\mathcal A}
\newcommand{\M}{\mathcal M}

\renewcommand{\D}{\mathcal D}
\newcommand{\Arf}{\mathit{Arf}}
\newcommand{\Proj}{\mathrm{Proj}}
\newcommand{\End}{\mathrm{End}}
\newcommand{\Hom}{\mathrm{Hom}}
\newcommand{\Mat}{\mathrm{M}}

\begin{document}

\begin{asciiabstract}
Bing doubling is an operation which produces a 2-component boundary
link B(K) from a knot K.  If K is slice, then B(K) is easily
seen to be boundary slice. In this paper, we investigate whether the
converse holds. Our main result is that if B(K) is boundary slice,
then K is algebraically slice.  We also show that the Rasmussen
invariant can tell that certain Bing doubles are not smoothly slice.
\end{asciiabstract}

\begin{htmlabstract}
Bing doubling is an operation which produces a 2-component boundary
link B(K) from a knot K.  If K is slice, then B(K) is easily
seen to be boundary slice. In this paper, we investigate whether the
converse holds. Our main result is that if B(K) is boundary slice,
then K is algebraically slice.  We also show that the Rasmussen
invariant can tell that certain Bing doubles are not smoothly slice.
\end{htmlabstract}

\begin{abstract}
Bing doubling is an operation which produces a 2--component boundary
link $B(K)$ from a knot $K$.  If $K$ is slice, then $B(K)$ is easily
seen to be boundary slice. In this paper, we investigate whether the
converse holds. Our main result is that if $B(K)$ is boundary slice,
then $K$ is algebraically slice.  We also show that the Rasmussen
invariant can tell that certain Bing doubles are not smoothly slice.
\end{abstract}

\maketitle

\section*{Introduction}

Bing doubling \cite{Bing} is a standard construction which, given a knot $K$ in $S^3$, produces a 2--component
oriented link $B(K)$ as illustrated below. (See \fullref{section:Bing} for a precise definition.)

\begin{figure}[htbp]
\labellist\small\hair 2.5pt
\pinlabel {$K$} at -30 210
\pinlabel {$B(K)$} at 350 210
\endlabellist
\centerline{\psfig{file=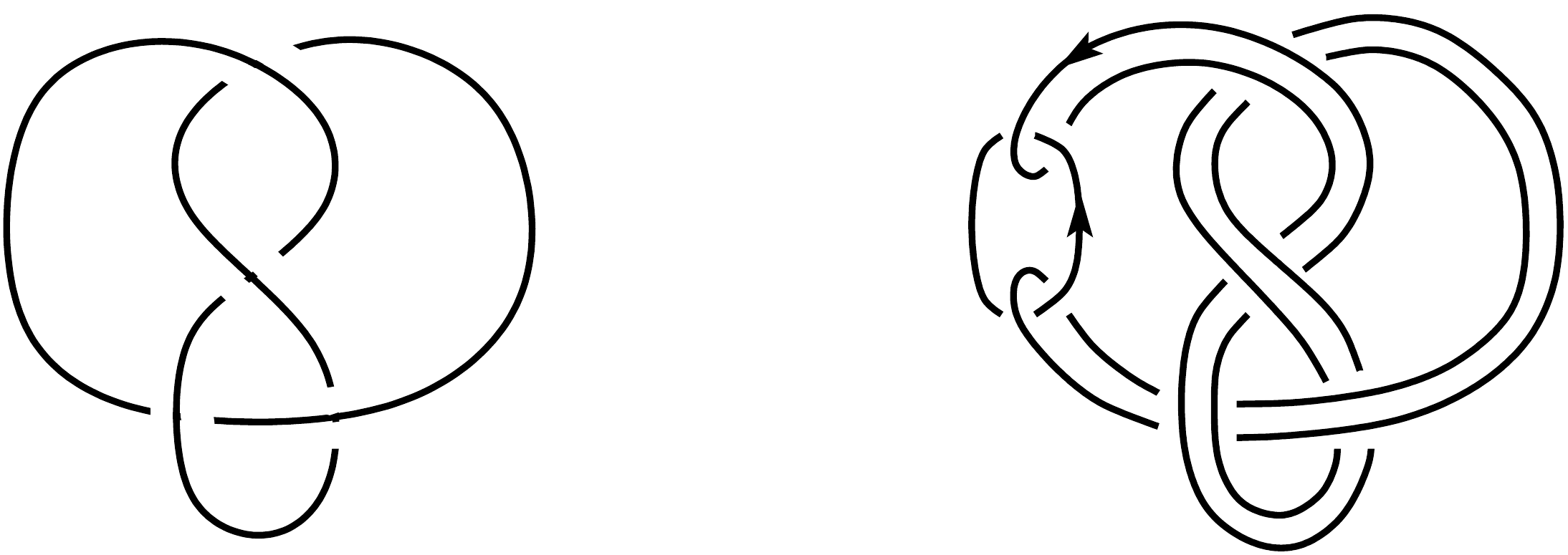,height=2.8cm}}
\caption{The figure eight knot and its Bing double}
\label{fig:ex}
\end{figure}

\noindent It is easy to check that if the knot $K$ is slice, then the link $B(K)$ is slice.
Does the converse hold~? An affirmative answer to this question seems out of reach.
However, a result obtained independently by S Harvey \cite{Har} and P Teichner
\cite{Tei} provides a first step in this direction.

Recall that the Levine--Tristram signature of a knot $K$ is the function $\sigma_K\co S^1\to\Z$ defined as follows:
for $\omega\in S^1$, $\sigma_K(\omega)$ is given by the signature of the Hermitian matrix $(1-\omega)A+(1-\overline\omega)A^*$,
where $A$ is a Seifert matrix for $K$ and $A^*$ denotes the transposed matrix.

\begin{thmintro}[Harvey \cite{Har}, Teichner \cite{Tei}]
If $B(K)$ is slice, then the integral over $S^1$ of the Levine--Tristram signature of $K$ is zero.
\end{thmintro}

For example, this shows that the Bing double of the trefoil knot is not slice. However, it says nothing about the Bing double of the figure eight knot
(\fullref{fig:ex}), as $\sigma_K$ is identically zero. Let us also
mention a conjecture of R Schneiderman and P Teichner: if
$B(K)$ is slice, then the Arf invariant of $K$ is zero. They hope to prove this statement using tree-valued intersections of Whitney
towers.

In this paper, we obtain a couple of results in the same direction.
It is known that Bing doubles are boundary links: their components bound disjoint Seifert surfaces (see also \fullref{prop:boundary}).
Furthermore, if a knot is slice, then its Bing double is not only slice, but boundary slice: the two slicing discs bound disjoint
oriented 3--manifolds in the 4--ball (see also \fullref{prop:conc}). Therefore, it makes sense to consider the following variation of the original
question: is $B(K)$ boundary slice only when $K$ is slice~?
Note that there is no known example of a boundary link which is slice and not boundary slice. Hence, it may turn out that both questions are in fact equivalent.

Using D Sheiham's work on boundary link concordance \cite{She}, we obtain the following result.

\begin{thmintro}
If $B(K)$ is boundary slice, then $K$ is algebraically slice.
\end{thmintro}

This shows in particular that the Bing double of the figure eight knot is not boundary slice, a fact known to experts like K Ko \cite{Komail}.

In a final section, we show that the Rasmussen concordance invariant \cite{Ras} (more specifically, its generalization to links
due to A Beliakova and S Wehrli \cite{B-W})
can be used to study the sliceness of Bing doubles. Indeed, we prove that if the Thurston--Bennequin invariant of a knot $K$ is non-negative, then $s(B(K))=1$.
Since the Rasmussen invariant of the 2--component unlink is $-1$, it follows that such a Bing double is not smoothly slice. This latter statement can 
also be obtained from results of L Rudolph \cite{Ru2}. Nevertheless, our computations show that this fact follows
very easily from elementary properties of the Rasmussen invariant. Note that there are examples of topologically slice knots with non-negative
Thurston--Bennequin invariant (see \fullref{section:s}). Bing doubling these knots gives the following result.
 
\begin{thmintro}
There exist links that are topologically slice but not smoothly slice, and whose components are trivial.
\end{thmintro}

\medskip

The paper is organized as follows. In \fullref{section:Bing}, we give a precise definition of Bing doubling and recall some of its
well-known properties. We also check that many `classical' concordance invariants are useless for studying Bing doubles.
(This is an attempt to convince the reader that the
question under study is non-trivial). \fullref{section:boundary} starts with a brief survey of boundary link concordance theory, including Sheiham's
\cite{She}. We then prove our main result. In \fullref{section:s}, we recall several properties of the Rasmussen invariant, and show that it can be
used to detect some Bing doubles that are not smoothly slice.

\section{The Bing double of a knot, first properties}\label{section:Bing}

Let $K$ be an oriented knot in $S^3$, and let $\N(K)$ denote a closed tubular neighborhood of $K$ in $S^3$. Recall that a pair $m,\ell$ of
oriented simple closed curves in $\partial \N(K)$ is a {\em standard meridian and longitude\/} for $K$ if $[m]=0$, $[\ell]=[K]$ in $H_1(\N(K))$,
and $lk(m,K)=1$, $lk(\ell,K)=0$, where $lk(\cdot,\cdot)$ denotes the linking number. Note that such a pair is unique up to isotopy.

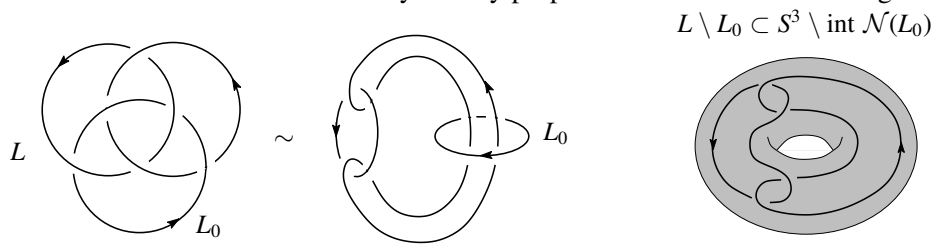
\begin{figure}[htbp]
\labellist\small\hair 2.5pt
\pinlabel {$L$} at -40 160
\pinlabel {$L_0$} at 290 30
\pinlabel {$L_0$} at 890 190
\pinlabel {$\sim$} at 420 180
\pinlabel {$L\setminus L_0\subset S^3\setminus \int\,\N(L_0)$} at 1320 380
\endlabellist
\centerline{\psfig{file=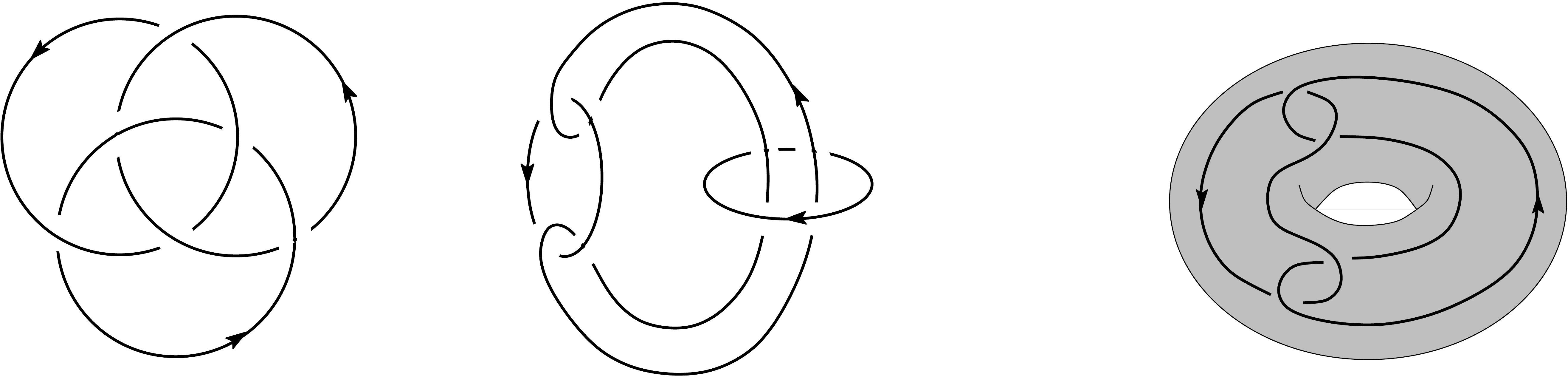,height=2.8cm}}
\caption{The Borromean rings $L$, and the solid torus $S^3\setminus \int\,\N(L_0)$}
\label{fig:Borromean}
\end{figure}

Let $L$ denote the Borromean rings, arbitrarily oriented, and let $m_0,\ell_0\subset \partial \N(L_0)$
be a standard meridian and longitude for some component $L_0$ of $L$. The {\em Bing double of $K$\/} is the 2--component oriented link $B(K)$ given by
the image of $L\setminus L_0$ in $(S^3\setminus \int\,\N(L_0))\cup(S^3\setminus \int\,\N(K))=S^3$, where the pasting homomorphism maps $m$ onto $\ell_0$ and
$\ell$ onto $m_0$. \fullref{fig:Borromean} should convince the reader that this definition corresponds to the illustration given in \fullref{fig:ex}.
Note that the isotopy type of the oriented link $B(K)$ depends neither on the choice of the component $L_0$ of $L$,
nor on the orientation of $L$ and $K$. This follows from obvious symmetry properties of the Borromean rings.

Recall that an oriented link is a {\em boundary link\/} if its components bound disjoint Seifert surfaces.
The following seems well known.

\begin{prop}\label{prop:boundary}
Bing doubles are boundary links.
\end{prop}
\begin{proof}
Consider the pair $L_1\cup L_2:=L\setminus L_0\subset S^3\setminus \int\,\N(L_0)=S^1\times D^2$ illustrated in \fullref{fig:Borromean},
where the solid torus is parametrized using the standard meridian and longitude for $L_0$. Note that there are two disjoint genus $0$
surfaces $P_1$ and $P_2$ in $S^1\times D^2$ such that the boundary of $P_i$ consists of $L_i$ together with two longitudes of the solid torus,
for $i=1,2$.  This is illustrated below.

\begin{figure}[htbp]
\centerline{\psfig{file=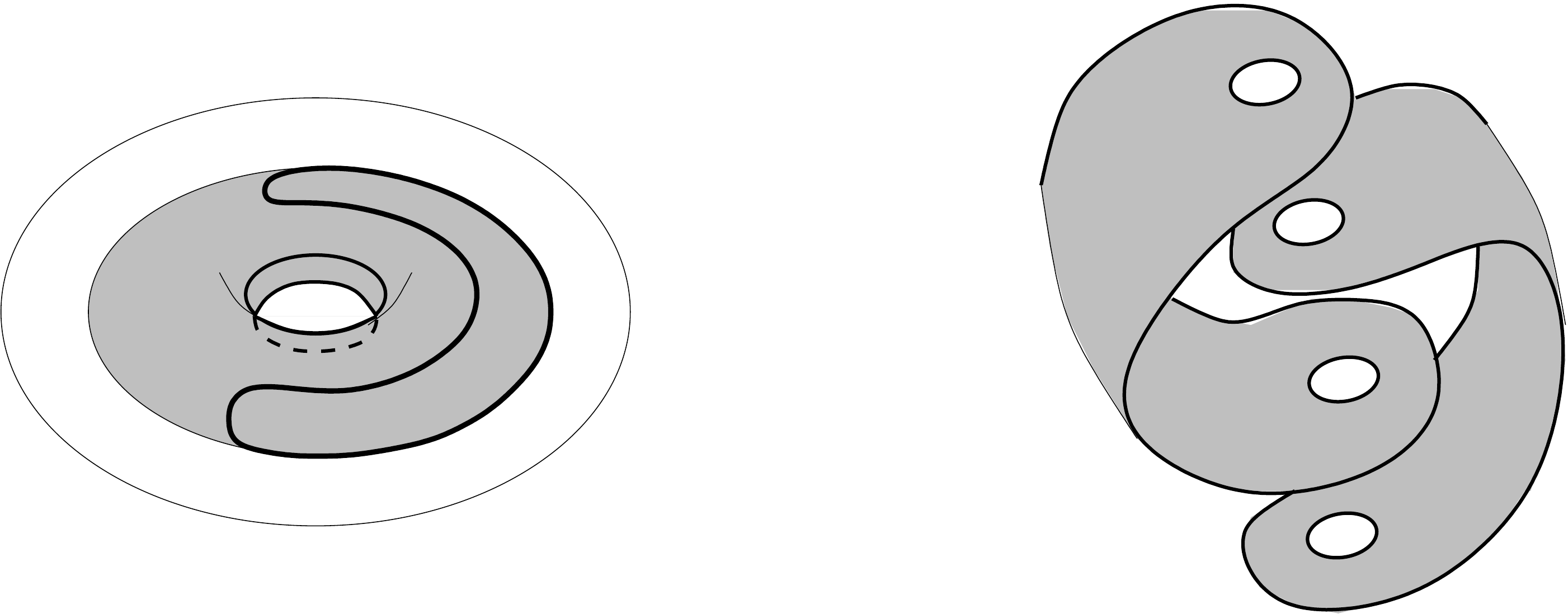,height=3cm}}
\caption{On the left, one of the $P_i$'s in the solid torus. On the right, an illustration of how $P_1$ and $P_2$ are embedded with
respect to each other.}
\label{fig:boundary}
\end{figure}

\noindent The pasting homeomorphism $h\co\partial(S^1\times D^2)\to\partial\N(K)$ maps the longitudes of $S^1\times D^2$ onto standard longitudes of
$\N(K)$, that is, parallel unlinked copies of $K$. These parallel copies bound disjoint Seifert surfaces (parallel copies of a fixed Seifert surface
for $K$). Pasting $h(P_1\sqcup P_2)$ with $4$ parallel Seifert surfaces for $K$, we obtain two disjoint Seifert surfaces for the components of $B(K)$.
\end{proof}

Of course, if $K$ and $K'$ are isotopic knots, then $B(K)$ and $B(K')$ are isotopic oriented links. Is the converse also true~? The answer is yes.
Indeed, the Jaco--Shalen--Johannson decomposition theorem implies that two knots $K$ and $K'$ are isotopic if and only if $B(K)$ and $B(K')$ are.
We refer to Hatcher \cite{Hat} for a beautiful exposition of the
JSJ--decomposition theorem, and to Budney \cite{Bud} for a survey of its consequences for knots and
links in $S^3$. The fact mentioned above can be understood as a special case of \cite[Proposition 4.31]{Bud}.

The aim of this paper is to address a 4--dimensional analogue of the question above. Recall that two $m$--component links $L_0$ and $L_1$ in $S^3$ are
{\em concordant\/} if there is a proper oriented locally flat submanifold $C\subset S^3\times[0,1]$, homeomorphic to $m$ copies of
$S^1\times[0,1]$, such that $C\cap S^3\times t=L_t$ for $t=0,1$. If these are boundary links and the concordance together with some disjoint Seifert
surfaces on top and bottom bound $m$ disjoint oriented 3--manifolds in $S^3\times[0,1]$, then the links are {\em boundary concordant\/}.
A (boundary) link is {\em (boundary) slice\/} if it is (boundary) concordant to the trivial link.
Note that a knot is slice if and only if it is boundary slice. (Indeed, if $V$ is a Seifert surface for a knot that has a
slicing disk $D$ in the 4--ball $B^4$, then elementary obstruction theory shows that $V\cup D$ bounds an oriented 3--manifold in $B^4$.)
Note also that there is no known example of a boundary link which is slice and not boundary slice.

\begin{prop}\label{prop:conc}
If two knots $K$ and $K'$ are concordant, then the links $B(K)$ and $B(K')$ are boundary concordant.
\end{prop}
\begin{proof}
Let $K$ and $K'$ be concordant knots with Seifert surfaces $V$ and $V'$. Fix a concordance
$C\co S^1\times[0,1]\hookrightarrow S^3\times[0,1]$ between $K$ and $K'$, and let $W$ be an oriented
3--manifold in $S^3\times[0,1]$ such that $\partial W=V\cup C\cup -V'$. Let $\tilde C\co S^1\times D^2\times[0,1]\hookrightarrow S^3\times[0,1]$
be a parametrization of a tubular neighborhood of $C$, such that $(*\times\partial D^2,S^1\times*)\times t$ maps to a standard meridian and longitude
for $K$ if $t=0$ and for $K'$ if $t=1$. Finally, let $P$ be a pair of pants, and let $\varphi\co P\sqcup P\hookrightarrow S^1\times D^2$ be a parametrization of the
embedding illustrated in \fullref{fig:boundary}. Consider the oriented 3--manifold $\widetilde W_1\sqcup\widetilde W_2$ given by the
image of the embedding
\[
(P\sqcup P)\times[0,1]\stackrel{\varphi\times id}{\longrightarrow}S^1\times D^2\times[0,1]\stackrel{\tilde C}{\longrightarrow}S^3\times[0,1].
\]
Note that for $i=1,2$, $\partial\widetilde W_i$ consists of a concordance between the $i^{\mathrm{th}}$ component of $B(K)$ and $B(K')$,
together with two parallel copies of $C$. Let $W_i$ denote the 3--manifold obtained by pasting $\widetilde W_i$ with
two parallel copies of $W$. The disjoint 3--manifolds $W_1$, $W_2$ provide a boundary concordance between $B(K)$ and $B(K')$.
\end{proof}
\begin{cor}\label{cor:slice}
If a knot $K$ is slice, then the link $B(K)$ is boundary slice.
\end{cor}

This latter result is well-known to experts. It motivates the following:

\begin{qu}
If $B(K)$ is a slice link, is $K$ necessarily a slice knot~?
\end{qu}

To illustrate the difficulty of this problem, let us go through a list of obstructions to the sliceness of links. We shall see
that all these obstructions vanish for {\em all\/} Bing doubles.

\medskip

\noindent{\bf The (multivariable) Alexander polynomial}\qua
It is well-known that the Alexander polynomial of a slice knot $K$ is of the form $\Delta_K(t)=f(t)f(t^{-1})$ for some $f(t)\in\Z[t]$.
A Kawauchi \cite{Kaw} generalized this result in the following way. Given an $m$--component ordered link $L$, let $\A(L)$ denote its Alexander module over
$\Z[t_1^{\pm 1},\dots,t_m^{\pm 1}]$, and let $\widetilde\Delta_L$ be the greatest common divisor of the order
ideal of the torsion of $\A(L)$. If $L$ is slice, then $\widetilde\Delta_L(t_1,\dots,t_m)=f(t_1,\dots,t_m)f(t^{-1}_1,\dots,t^{-1}_m)$
for some $f(t_1,\dots,t_m)\in\Z[t_1,\dots,t_m]$.

Let us now describe briefly how to compute the Alexander module of a
Bing double, referring to Cooper \cite{Coo} and Cimasoni--Florens \cite{C-F} for details.
Recall that a {\em C--complex\/} for a 2--component link $L=L_1\cup L_2$
is given by two Seifert surfaces $V_1$ for $L_1$ and $V_2$ for $L_2$ that intersect only along clasps. Given such a C--complex $V=V_1\cup V_2$,
one can define two `generalized Seifert forms' on $H_1(V)$, giving two Seifert matrices, say $A$ and $A'$. A presentation matrix for $\A(L)$
can be obtained by some $\Z[t_1^{\pm 1},t_2^{\pm 1}]$--linear combination of these matrices and their transpose.
This computational method is very efficient in the case of a Bing double $L=B(K)$. Indeed, there is an
obvious C--complex $V$ for $B(K)$ given by two discs intersecting along two clasps. Then, $H_1(V)=\Z$, and since the Bing double is untwisted, we get
$A=A'=\begin{pmatrix}0\end{pmatrix}$. Therefore, $\A(B(K))=\Z[t_1^{\pm 1},t_2^{\pm 1}]$ for any knot $K$, so $B(K)$ has the Alexander module
of a trivial link. In particular, $\widetilde\Delta_{B(K)}(t_1,t_2)=1$ for all knots $K$.

\medskip

\noindent{\bf The (multivariable) Levine--Tristram signature}\qua
The Levine--Tristram signature of a knot $K$ is the function $\sigma_K\co S^1\to\Z$ defined as follows: for $\omega\in S^1$, $\sigma_K(\omega)$ is given
by the signature of the Hermitian matrix $(1-\omega)A+(1-\overline\omega)A^*$, where $A$ is a Seifert matrix for $K$. If $K$ is a
slice knot, then $\sigma_K(\omega)=0$ for all $\omega\in S^1$ such that $\Delta_K(\omega)\neq 0$. This invariant admits a multivariable generalization:
for a 2--component ordered link $L$, it consists of a function $\sigma_L\co S^1\times S^1\to\Z$. The integer $\sigma_L(\omega_1,\omega_2)$ is the
signature of a Hermitian matrix given by some $\Cc$--linear combination (depending on $\omega_1$ and $\omega_2$) of the Seifert matrices $A$ and $A'$
for $L$. If $L$ is a slice link, then $\sigma_L(\omega_1,\omega_2)=0$ for all $\omega_1,\omega_2\in S^1$ such that
$\widetilde\Delta_L(\omega_1,\omega_2)\neq 0$. (Here again, we refer to \cite{Coo} for the case of 2--component links, and to \cite{C-F} for the
general case.)

As mentioned in the previous paragraph, there is a choice of C--complex for $L=B(K)$ such that $A=A'=\begin{pmatrix}0\end{pmatrix}$.
Therefore, $\sigma_{B(K)}=0$ for all knots $K$, so this invariant does not tell us anything about Bing doubles.

\begin{rem}
In the construction of the Bing double of $K$, if one replaces the standard longitude for $K$ by some longitude that links $K$ (say, $t$ times), then the
result is the $t$--twisted Bing double $B(K,t)$. One easily checks that $\sigma_{B(K,t)}(\omega_1,\omega_2)$ is equal to the sign of $t$ for all
$\omega_1,\omega_2\neq 1$. Hence, the $t$--twisted Bing double of a knot is {\em never\/} slice if $t\neq 0$ (unlike some $t$--twisted Whitehead doubles).
This is the reason why we restrict ourselves to the study of untwisted Bing doubles.
\end{rem}

\medskip

\noindent{\bf The Arf invariant}\qua
Let $\alpha_2\co H_1(V;\Z_2)\times H_1(V;\Z_2)\to\Z_2$ denote the Seifert form of a knot $K$ reduced modulo 2. The Arf invariant $\Arf(K)$
of $K$ is defined as the Arf invariant of the non-singular quadratic form $q\co H_1(V;\Z_2)\to\Z_2$ given by
$q(x)=\alpha_2(x,x)$. If $K$ is a slice knot, then $\Arf(K)=0$. If $L$ is an oriented link with components $\{L_i\}$ that satisfy
$lk(L_i,L\setminus L_i)\equiv 0\pmod{2}$ for all $i$, then its reduced Seifert form $\alpha_2$ induces a well-defined non-singular quadratic form on
$H_1(V;\Z_2)/i_\ast H_1(\partial V;\Z_2)$. The Arf invariant of $L$ is then defined as the Arf invariant of this quadratic form.

The Bing double of a knot obviously satisfies the condition above and not surprisingly, $\Arf(B(K))=0$ for all knots $K$.
There are numerous ways to check this fact. For example, K Murasugi \cite{Mur} showed that if $L=L_1\cup L_2$ satisfies
$lk(L_1,L_2)\equiv 0\pmod{2}$, then
\[
\Arf(L)\equiv \Arf(L_1)+\Arf(L_2)+\frac{1}{2}\frac{d^2}{dt^2}\Delta_L(t,t)\Big\vert_{t=1}\pmod{2}.
\]
In the case of $L=B(K)$, both components are trivial and $\Delta_L(t_1,t_2)=0$, leading to the result.

\medskip

\noindent{\bf The Milnor $\overline\mu$--invariants}\qua
As we saw in \fullref{prop:boundary}, $B(K)$ is a boundary link.
This implies the vanishing of another obstruction to the sliceless of $B(K)$: Milnor's $\overline\mu$--invariants.

\medskip

\noindent{\bf The slice Bennequin inequality}\qua
Up to now, we have been concerned with the {\em topological\/} concordance of links. If one requires the concordance to be {\em smooth\/}, one
gets a stronger equivalence relation: there are links which are topologically slice (bound disjoint locally flat discs in $B^4$) but not
smoothly slice (do not bound disjoint smooth discs in $B^4$). We shall now investigate an obstruction to the smooth sliceness of links.
Given an oriented link $L$, let $\chi_s(L)$ denote the greatest Euler characteristic $\chi(F)$ of an oriented surface $F$ (with no closed component)
smoothly embedded in $B^4$ with boundary $L$. The {\em slice Bennequin
  inequality\/} (Rudolph \cite{Ru1}) asserts that for every braid $\beta\in B_n$,
\[
\chi_s(\what\beta)\le n-\omega(\beta),
\]
where $\what\beta$ denotes the closure of $\beta$ and $\omega(\beta)$ its writhe (ie, the number of positive crossings minus the number of negative ones).

This inequality does not tell us anything about Bing doubles. Indeed, let $L=B(K)=L_1\cup L_2$ and fix $\beta\in B_n$ such that $\what\beta=L$.
Write the geometric braid $\beta$ as a union $\beta=\beta_1\cup\beta_2$ with $\what\beta_i=L_i$
and $\beta_i\in B_{n_i}$ for $i=1,2$. Clearly, $n_1+n_2=n$ and $\omega(\beta)=\omega(\beta_1)+\omega(\beta_2)+2\,lk(L_1,L_2)=\omega(\beta_1)+\omega(\beta_2)$.
Since each $L_i$ is the unknot, the slice Bennequin inequality applied to $\beta_i$ gives $1\le n_i-\omega(\beta_i)$. Summing over $i=1,2$ gives
$2\le n-\omega(\beta)$. This means that, in the best case, the slice Bennequin inequality for $\beta$ will read $\chi_s(L)\le 2$. This does not say anything about $L$,
as this inequality holds for any 2--component link.

\medskip

To conclude this panorama, it should be mentioned that invariants of Cha--Ko \cite{C-K} and Friedl \cite{Fri} do detect the Bing double of some
knots with non-trivial Levine--Tristram signature. This is also true
for the $L^2$--signatures of Harvey \cite{Har}, as mentioned in the introduction.
However, none of these invariants can detect the Bing double of a torsion element in the knot concordance group.

\section{Boundary sliceness of Bing doubles}\label{section:boundary}

Since Bing doubles are boundary links, and since $B(K)$ is boundary slice whenever $K$ is slice,
it makes sense to consider the following variation of our original problem:

\begin{qu}
If $B(K)$ is boundary slice, is $K$ necessarily slice~?
\end{qu}

Let us recall once again that there is no known example of a boundary link which is slice and not boundary slice.
Therefore, both questions might turn out to be equivalent. Nevertheless, we shall be more successful with this version.

\medskip

In order to state our results, let us recall several standard facts
about boundary link concordance. We refer to Ko \cite{Ko} for proofs and further details.

Let $B(m)$ denote the set of boundary concordance classes of $m$--component boundary links in $S^3$.
A pair $(L,V)$ consisting of a boundary link $L=L_1\cup\dots\cup L_m$ together with disjoint Seifert surfaces $V=V_1\cup\dots\cup V_m$ is called a
boundary pair. Two boundary pairs $(L,V)$ and $(L',V')$ are said to be concordant if there are disjoint oriented 3--manifolds $W_1,\dots,W_m$
in $B^4$ such that $\partial W_i=V_i\cup C_i\cup -V_i'$, where $C_i$ is a concordance between $L_i$ and $L_i'$. The set of concordance classes of
boundary pairs is denoted by $C(B_m)$. Then, the obvious surjective map $C(B_m)\twoheadrightarrow B(m)$ turns out to be a bijection for $m\le 2$.
(This is false for $m>2$.)

Fix an $m$--component boundary pair $(L,V)$. For $i,j=1,\dots,m$, consider the bilinear
pairing $H_1(V_i)\times H_1(V_j)\to\Z$ defined by $(x,y)\mapsto lk(x^+,y)$, where $x^+$ denotes the 1--cycle $x$ pushed in the positive normal direction
of $V_i$. Let $A_{ij}$ denote the matrix of this pairing with respect to some fixed bases of $H_1(V_i)$ and $H_1(V_j)$. Note that $A_{ij}=A_{ji}^*$ if $i\neq j$ and 
$A_{ii}-A_{ii}^*$ is nothing but the intersection matrix of $V_i$. Furthermore, if the boundary pair $(L,V)$ is slice, then for suitable bases of the
$H_1(V_i)$'s, each $A_{ij}$ is metabolic: the upper left quadrant is zero. In this case, the collection $A=\{A_{ij}\}$ (and the boundary pair $(L,V)$)
are said to be {\em algebraically slice\/}. 
This motivates the following definition: consider the set of collections $A=\{A_{ij}\}_{i,j=1}^m$ of integral matrices such that
$A_{ij}=A_{ji}^*$ if $i\neq j$ and $A_{ii}-A_{ii}^*$ is unimodular. Let us say that two such collections $A$ and $B$ are equivalent if the
{\em block sum\/} $A\oplus(-B)$ defined by $(A\oplus(-B))_{ij}=A_{ij}\oplus(-B_{ij})$ is algebraically slice. Then, the set of equivalence classes
forms an abelian group $G^m(\Z)$ under the block sum.

The Seifert matrix construction described above defines a map $\psi_m\co C(B_m)\to G^m(\Z)$ which is onto but not injective.
Recall that there are canonical bijections $C(B_1)\cong B(1)\cong\mathcal{C}$ and $C(B_2)\cong B(2)$, where $\mathcal{C}$
is the classical knot concordance group. Therefore, we have an epimorphism
$\psi_1\co\mathcal{C}\twoheadrightarrow G^1(\Z)$ and a surjective map $\psi_2\co B(2)\twoheadrightarrow G^2(\Z)$.

\medskip

We are finally ready to state our main result.

\begin{thm}\label{thm:boundary}
{\rm(i)}\qua The Bing doubling operation induces a map $\mathcal{C}\to B(2)$ and a homomorphism $\varphi\co G^1(\Z)\to G^2(\Z)$
such that the following diagram commutes:
\[
\xymatrix{
\mathcal{C}\ar@{->>}[d]_{\psi_1}\ar[r]& B(2)\ar@{->>}[d]^{\psi_2}\\
G^1(\Z)\ar[r]^{\varphi}& G^2(\Z)
}
\]
{\rm(ii)}\qua The homomorphism $\varphi$ is injective.
\end{thm}

\begin{cor}
If $B(K)$ is boundary slice, then $K$ is algebraically slice.
\end{cor}

\begin{proof}[Proof of $\mathbf{(i)}$]
The fact that Bing doubling induces a map $\mathcal{C}\to B(2)$ follows from \fullref{prop:conc}.
Let us define the homomorphism $\varphi\co G^1(\Z)\to G^2(\Z)$. Given a matrix $A$ such that $A-A^*$ is unimodular, set
\[
B(A)=\begin{pmatrix}
A_{11}&A_{12}\cr A_{21}&A_{22}
\end{pmatrix}
\]
where:
\begin{equation}\label{B(A)}
A_{11}=A_{22}=\begin{pmatrix}
A\phantom{^*}&A\phantom{^*}\cr A^*&A^*
\end{pmatrix}\quad\hbox{and}\quad
A_{12}=A^*_{21}=\begin{pmatrix}
A\phantom{^*}&A\cr A^*&A
\end{pmatrix}
\end{equation}
Note that $A_{12}=A^*_{21}$ and that $A_{11}-A^*_{11}=A_{22}-A^*_{22}=\begin{pmatrix}A-A^*&0\cr 0&A-A^*\end{pmatrix}$ is unimodular since $A-A^*$ is.
Hence, $B(A)$ defines an element of the group $G^2(\Z)$.
Let us now assume that $A$ is algebraically slice, say of size $2g$. This means that there exists some unimodular matrix $Q$ of size $2g$ such that
$QAQ^*$ is metabolic. Consider the unimodular matrix
\[
\widehat Q=\begin{pmatrix}
I_g&0&0&0\cr
0&0&I_g&0\cr
0&I_g&0&0\cr
0&0&0&I_g
\end{pmatrix}
\begin{pmatrix}
Q&0\cr 0&Q
\end{pmatrix}
\]
where $I_{g}$ denotes the identity matrix of size $g$. A straightforward computation shows that $\widehat QA_{ij}\widehat Q^*$ is metabolic for all $i,j$.
Hence, $B(A)$ is algebraically slice. Now, let $A$ and $A'$ be two Seifert matrices of respective size $2g$ and $2g'$. Consider the unimodular matrix:
\[
R=\begin{pmatrix}
I_{2g}&0&0&0\cr
0&0&I_{2g}&0\cr
0&I_{2g'}&0&0\cr
0&0&0&I_{2g'}
\end{pmatrix}
\]
One checks that
\[
\begin{pmatrix}R&0\cr 0&R\end{pmatrix}
B(A\oplus A')
\begin{pmatrix}R^*&0\cr 0&R^*\end{pmatrix}=B(A)\oplus B(A').
\]
These facts imply that $A\mapsto B(A)$ defines a homomorphism $\varphi\co G^1(\Z)\to G^2(\Z)$.
It remains to check that the diagram commutes. Let $(K,V)$ be a 1--component boundary
pair (that is, a knot with a Seifert surface), and let $A$ denote the associated Seifert matrix with respect to some basis $b$ of $H_1(V)$.
Let $B(V)=B(V)_1\sqcup B(V)_2$ denote the boundary Seifert surface for $B(K)$ constructed in the proof of \fullref{prop:boundary}.
Clearly, $H_1(B(V)_1)=H_1(B(V)_2)=H_1(V)\oplus H_1(V)$. Consider the right hand side of \fullref{fig:boundary}: the four glued in
Seifert surfaces are parallel, but their orientations alternate. It follows that the Seifert matrix with respect to the basis $b\cup b$ of $H_1(B(V)_1)$ and $H_1(B(V)_2)$ is given by $B(A)$.
\end{proof}

Our proof the injectivity of $\varphi$ makes extensive use of D Sheiham's interpretation of $G^m(\Z)$ as the Witt group of
the representation category of some ring $P_m$. We shall now quickly review the notions and results involved, referring to \cite{She} for
further details and proofs.

Fix a commutative ring $\Lambda$ and a ring $R$ with an involution $r\mapsto\overline{r}$.
Let $(R-\Lambda)-\Proj$ denote the category of representations $\rho\co R\to\End_\Lambda(M)$, where $M$ is a
finitely generated projective $\Lambda$--module. This is a Hermitian category via $(M,\rho)\mapsto (M^*,\rho^*)$, where $M^*=\Hom_\Lambda(M,\Lambda)$
and $\rho^*\co R\to\End_\Lambda(M^*)$ is given by $\rho^*(r)(\xi)\co x\mapsto\xi(\rho(\overline{r})(x))$ for $r\in R$, $\xi\in M^*$
and $x\in M$.

Fix a sign $\epsilon=\pm 1$. An {\em $\epsilon$--hermitian form\/} $(M,\phi)$ is a morphism $\phi\co M\to M^*$ in $(R-\Lambda)-\Proj$ such that
$\phi^*=\epsilon\phi$. (Here, $(M^*)^*$ is identified with $M$.) An object $M$ is {\em $\epsilon$--self dual\/} if there exists an $\epsilon$--hermitian form $(M,\phi)$ which is non-singular,
ie an isomorphism. Finally, recall that an $\epsilon$--hermitian form $(M,\phi)$ is called {\em metabolic\/} if there is a direct summand
$j\co L\hookrightarrow M$ such that $L=\ker(j^*\phi\co M\to L^*)$. The {\em Witt group\/} of the Hermitian category
$(R-\Lambda)-\Proj$ is the abelian group $W^\epsilon(R-\Lambda)$ with one generator $[M,\phi]$ for each isomorphism class of non-singular $\epsilon$--hermitian
form $(M,\phi)$ subject to the relations
\[
\begin{cases}
[M',\phi']=[M,\phi]+[M'',\phi'']& \text{if $(M',\phi')\cong(M,\phi)\oplus(M'',\phi'')$;}\\
[M,\phi]=0 & \text{if $(M,\phi)$ is metabolic.}
\end{cases}
\]
Note that $W^\epsilon(\Z-\Lambda)$ is simply the Witt group $W^\epsilon(\Lambda)$ of the ring $\Lambda$.

Let us now focus on some particular choice of the ring $R$. Let $m\ge 1$ be an integer, and let $P_m$ be the ring
\[
P_m=\Z\left<s,\pi_1,\dots,\pi_m\,\big|\,\text{$\pi_i^2=\pi_i$, $\pi_i\pi_j=0$ for $i\neq j$, $\textstyle{\sum_{i=1}^m\pi_i=1}$}\right>
\]
endowed with the involution induced by $s\mapsto\overline{s}=1-s$ and $\pi_i\mapsto\overline{\pi_i}=\pi_i$ for $i=1,\dots,m$.
Given any commutative ring $\Lambda$, Sheiham's idea is to interpret the group $G^m(\Lambda)$ as the Witt group $W^{-1}(P_m-\Lambda)$.
Indeed, an element in $G^m(\Lambda)$ given by a collection of matrices $A=\{A_{ij}\}_{i,j=1}^m$ can be understood
as the equivalence class of an $(m+2)$--tuple $(M,\pi_1,\dots,\pi_m,\lambda)$, where $M$ is a finitely generated projective $\Lambda$--module,
$\pi_1,\dots,\pi_m$ is a set of orthogonal idempotents in $\End_\Lambda(M)$ and $\lambda\co M\to M^*$ is an $\Lambda$--module homomorphism such
that $\lambda-\lambda^*$ is an isomorphism which commutes with each $\pi_i$.
($A_{ij}$ should be thought of as the matrix of a homomorphism $\pi_iM\to(\pi_jM)^*$, where the $\pi_iM$'s are free $\Lambda$--modules;
the collection of matrices $A$ then defines a homomorphism $\lambda\co M\to M^*$, where $M=\bigoplus_{i=1}^m\pi_iM$.)
Given such an $(m+2)$--tuple, set
$\phi=\lambda-\lambda^*$ and define $\rho\co P_m\to\End_\Lambda(M)$ by $\rho(s)=\phi^{-1}\lambda$, $\rho(\pi_i)=\pi_i$. One easily
checks that $(M,\phi)$ is a non-singular $(-1)$--hermitian form. The only non-trivial point is the equality
$\phi(sx)(y)=\phi(x)(\overline{s}y)$ for $x,y\in M$, where $sx$ stands for $\rho(s)(x)$. It follows from the
following equations:
\begin{align*}
\phi(sx)(y)+\phi(x)(sy)&=\phi\phi^{-1}\lambda(x)(y)+\phi(x)(\phi^{-1}\lambda(y))\\
	&=\lambda(x)(y)+\lambda^*(\phi^{-1})^*\phi(x)(y)=(\lambda-\lambda^*)(x)(y)=\phi(x)(y).
\end{align*}
It turns out that this construction induces an isomorphism
\[
G^m(\Lambda)\stackrel{\kappa}{\longrightarrow}W^{-1}(P_m-\Lambda),\quad[M,\pi_1,\dots,\pi_m,\lambda]\longmapsto[M,\phi]
\]
which is natural with respect to $\Lambda$ (see \cite[Lemma 3.31]{She}). 

The next step makes use of so-called {\em Hermitian devissage\/}. For the Hermitian category
$(R-k)-\Proj$, where $k$ is a field, it can be stated as follows. Let $\overline{\M}^s(R-k,\epsilon)$ denote the set of isomorphism classes of
simple $\epsilon$--self-dual objects in $(R-k)-\Proj$. Then, there is a canonical isomorphism of Witt groups
\[
W^\epsilon(R-k)\cong\bigoplus_{M\in\overline{\M}^s(R-k,\epsilon)}W^\epsilon((R-k)|_M),
\]
where $(R-k)|_M$ denotes the (Hermitian) full subcategory of $(R-k)-\Proj$ whose objects are isomorphic to 
a direct summand of $M^{\oplus d}$ for some $d$.

\bigskip
We conclude this brief exposition of Sheiham's work with the following special case of {\em Hermitian Morita equivalence\/}.
Fix a non-singular $\epsilon$--hermitian form $(M,b)$ in $(R-k)-\Proj$. Then, we have an isomorphism
\[
\begin{array}{ccl}
W^\epsilon((R-k)\vert_M)&\stackrel{\Theta_{M,b}}{\longrightarrow}&W^1(\End_{(R-k)}M)\cr
[N,\phi]&\longmapsto&[\Hom(M,N),\Phi_N\phi_*],
\end{array}
\]
where $\Hom(M,N)=\Hom_{(R-k)}(M,N)$ is endowed with the obvious structure of right $\End_{(R-k)}M$--module, and
$\Phi_N\phi_*\co\Hom(M,N)\to\Hom(M,N)^*$ is given by $(\Phi_N\phi_*)(\alpha)(\beta)=\epsilon b^{-1}\alpha^*\phi\beta$ for
$\alpha,\beta\in\Hom(M,N)$.
\bigskip

Finally, note that the natural homomorphism $G^m(\Z)\to G^m(\Q)$ is injective \cite[Lemma 1.11]{She}.
Summing up, we have the following sequence of homomorphisms:
\begin{align*}
G^m(\Z)&\hookrightarrow G^m(\Q)\stackrel{\kappa}{\longrightarrow}W^{-1}(P_m-\Q)\cong\bigoplus_{M}W^{-1}((P_m-\Q)|_M)\\
&\stackrel{p_M}{\twoheadrightarrow} W^{-1}((P_m-\Q)|_M)\stackrel{\Theta_{M,b}}{\longrightarrow}W^1(\End_{(P_m-\Q)}M),
\end{align*}
where the direct sum is over all $M\in\overline{\M}^s(P_m-\Q,-1)$, and $p_M$ denotes the canonical projection corresponding to $M$.

\bigskip

\noindent{\bf Proof of $\mathbf{(ii)}$}\qua
We shall now use these results to prove the second part of \fullref{thm:boundary}.
Recall that the Bing doubling map $\varphi\co G^1(\Z)\to G^2(\Z)$ (and similarly, the induced map $G^1(\Q)\to G^2(\Q)$) is given by
$[A]\mapsto[B(A)]=\begin{bmatrix}A_{11}&A_{12}\cr A_{21}&A_{22}\end{bmatrix}$, where the $A_{ij}$'s are given by \eqref{B(A)}.
This matrix $B(A)$ is block-congruent to (and therefore, represents the same class as) the matrix
\[
\widehat A:=\begin{pmatrix}
0&T&T&0\cr 0&A^*&T^*&A\cr T^*&T&0&T\cr 0&A^*&0&A^*
\end{pmatrix},
\]
where the letter $T$ stands for the unimodular matrix $A-A^*$. The isomorphism $\kappa\co G^1(\Q)\to W^{-1}(P_1-\Q)$ maps $[A]$ to $[M,\phi]$,
where $M$ is the $\Q$--vector space acted on by $A$, $\phi\co M\to M^*$ is given by the matrix $T=A-A^*$, and the action of $P_1=\Z[s]$
on $M$ is defined by the matrix $S=T^{-1}A$. Similarly, $\kappa\co G^2(\Q)\to W^{-1}(P_2-\Q)$ maps $[\widehat A]$ to
$[\widehat M,\what\phi]$, where $\widehat M=M^{\oplus 4}$ as a $\Q$--vector space, $\what\phi\co\widehat M\to(\widehat M)^*$ is given by the matrix
\begin{equation}\label{hatT}
\widehat T:=\widehat A-(\widehat A)^*=\begin{pmatrix}
0&T&0&0\cr T&T^*&0&0\cr 0&0&0&T\cr 0&0&T&T^*
\end{pmatrix}
\end{equation}
and the $P_2$--action on $\widehat M$ is defined as follows: for $x=(x_1,x_2,x_3,x_4)\in M^{\oplus 4}=\widehat M$, $\pi_1x=(x_1,x_2,0,0)$,
$\pi_2x=(0,0,x_3,x_4)$, and $sx=\widehat Sx$, where:
\begin{equation}\label{hatS}
\widehat S=(\widehat T)^{-1}\widehat A=\begin{pmatrix}
0&S&0&S\cr 0&I&I&0\cr -I&S&0&S\cr -I&I&0&I
\end{pmatrix}
\end{equation}

With these notations, we have the commutative diagram
\[
\xymatrix{
G^1(\Z)\ar@{^{(}->}[d]\ar[r]^\varphi& G^2(\Z)\ar@{^{(}->}[d]\\
G^1(\Q)\ar[r]\ar[d]^\cong_\kappa&G^2(\Q)\ar[d]_\kappa^\cong\\
W^{-1}(P_1-\Q)\ar[r]&W^{-1}(P_2-\Q)
}
\]
where the bottom homomorphism is given by $[M,\phi]\mapsto[\widehat M,\what\phi]$ as described above.
\bigskip

From now on, let us denote by $\C$ and $\D$ the Hermitian categories $(P_1-\Q)-\Proj$ and $(P_2-\Q)-\Proj$, respectively. 

\begin{lemma}\label{lemma:boundary1}
If $M$ is a simple object in $\C$, then $\widehat M$ is a simple object in $\D$.
\end{lemma}
\begin{proof}
First note that if $M$ is a simple representation over $\Q$ of $P_1=\Z[s]$, then $M=\Q[s]/(p)$ with $p\in\Q[s]$ some irreducible
polynomial in $\Q[s]$. If $S$ denotes a matrix of the action of $s\in\Z[s]$ on $M$ viewed as a $\Q$--vector space, then its characteristic polynomial $\det(S-sI)$ is equal to $p(s)$
up to multiplication by a non-zero scalar. Therefore, $\det(S-sI)$ must be irreducible in $\Q[s]$.

Let $n$ denote the dimension of $M$ over $\Q$, and let us assume that there is a $P_2$--submodule $N$ of $\widehat M$ with $N\neq 0$ and
$N\neq\widehat M$. In particular, the $\Q$--vector space $\widehat M=\pi_1\widehat M\oplus\pi_2\widehat M$ can be further decomposed into
$\pi_1\widehat M=\pi_1N\oplus N'_1$ and $\pi_2\widehat M=\pi_2N\oplus N'_2$. Set $d_i=\dim_\Q\pi_iN$ for $i=1,2$. Since $N$ is invariant
under the action of $s\in P_2$, there exists $P,Q\in GL_{2n}(\Q)$ such that
\begin{equation}\label{ABCD}
\begin{pmatrix}P&0\cr 0&Q\end{pmatrix}\widehat S\begin{pmatrix}P^{-1}&0\cr 0&Q^{-1}\end{pmatrix}=
\begin{pmatrix}A&0&C&0\cr *&*&*&*\cr D&0&B&0\cr *&*&*&*\end{pmatrix}
\end{equation}
for some matrices $A,B,C,D$ of respective size $d_1\times d_1$, $d_2\times d_2$, $d_1\times d_2$, and $d_2\times d_1$.
Here, $\widehat S$ denotes the matrix obtained from $S$ via \eqref{hatS}.

\medskip
{\bf Claim}\qua $d_1=d_2$

\noindent Since $M$ is simple, $p(s)=\det(S-sI)$ is irreducible in $\Q[s]$. In particular, $p(0)=\det(S)\neq 0$, so $S$ is invertible.
Therefore, $\begin{pmatrix}0&S\cr I&0\end{pmatrix}$ is invertible, as well as the matrix
$P\begin{pmatrix}0&S\cr I&0\end{pmatrix}Q^{-1}=\begin{pmatrix}C&0\cr *&*\end{pmatrix}$. Hence, the first $d_1$ lines of this matrix
are linearly independant. Since $C$ is of size $d_1\times d_2$, this implies that $d_1\le d_2$. Similarly, $0\neq p(1)=\det(S-I)$ implies
that $\begin{pmatrix}-I&S\cr -I&I\end{pmatrix}$ is invertible, as well as
$Q\begin{pmatrix}-I&S\cr -I&I\end{pmatrix}P^{-1}=\begin{pmatrix}D&0\cr *&*\end{pmatrix}$. Therefore, $d_2\le d_1$, proving the claim.
This also shows that the matrices $C$ and $D$ are invertible. Note that $d:=d_1=d_2$ satisfies $0<d<2n$, as $N\neq 0$ and
$N\neq\widehat M$.

Using this claim and the fact that $C$ is invertible, equation \eqref{ABCD} can be transformed into the following equalities:
\[
Q\begin{pmatrix}0&S\cr 0&I\end{pmatrix}=\begin{pmatrix}B&0\cr *&*\end{pmatrix}Q,\quad
Q\begin{pmatrix}I&0\cr I&0\end{pmatrix}=\begin{pmatrix}E&0\cr *&*\end{pmatrix}Q
\]
where $E=C^{-1}AC$. This implies that $Q=\begin{pmatrix}X&(E-I)X\cr *&*\end{pmatrix}$ for some matrix $X$ of size $d\times n$ which
satisfies the following equations:
\begin{equation}\label{BEX}
(B-I)(E-I)X=XS,\quad E(E-I)X=0,\quad BX=0.
\end{equation}

{\bf Claim}\qua $0<{\mathrm rank}\,X<n$

\noindent If $r:={\mathrm rank}\,X=0$, then $X=0$ so the first $d$ lines of $Q$ are zero.
This is impossible since $Q$ is invertible and $d>0$.
On the other hand, let us assume that $r=n$. This means that $X$ is the matrix of an injective linear map $M\hookrightarrow\Q^d$, which
we also denote by $X$. Since $S$ is invertible, the first equality in \eqref{BEX} implies that the rank of $(E-I)X$ is $n$. Consider
the subspace $V$ of $\Q^d$ consisting of the intersection of the images of $X$ and $(E-I)X$. Since both matrices have rank $n$ and since
$d<2n$, the space $V$ has positive dimension. Now, set
\[
W=\{w\in M\,\vert\,Xw=(E-I)Xm\hbox{ for some $m\in M$}\}. 
\]
Since $W=\{w\in M\,\vert\,Xw\in V\}$ and $X$ is injective, $\dim_\Q W=\dim_\Q V>0$. Given any $w\in W$, we have
$EXw=E(E-I)Xm=0$ by the second equation of \eqref{BEX}. Therefore,
\[
XSw=(B-I)(E-I)Xw=-BXw+Xw=Xw.
\]
Since $X$ has maximal rank, this implies that $Sw=w$ for all $w\in W\subset M$. As $W\neq 0$, this would imply that
$M$ is not a simple $\Q[s]$--module. Therefore, $r<n$, proving the claim.
\bigskip

Let us choose two matrices $T\in GL_d(\Q)$ and $U\in GL_n(\Q)$ such that $TXU=\begin{pmatrix}I_r&0\cr 0&0\end{pmatrix}=:\widetilde X$,
where $I_r$ denotes the identity matrix of size $r={\mathrm rank}\,X$. Equations \eqref{BEX} above give
\[
(\widetilde B-I)(\widetilde E-I)\widetilde X=\widetilde X\widetilde S,\quad \widetilde B\widetilde X=0,
\]
where $\widetilde E=TET^{-1}$, $\widetilde B=TBT^{-1}$ and $\widetilde S=U^{-1}SU$. Writing
$S=\begin{pmatrix}S_{11}&S_{12}\cr S_{21}&S_{22}\end{pmatrix}$ with $S_{11}$ and $S_{22}$ of respective size $r\times r$ and
$(n-r)\times(n-r)$, these equations imply that $S_{12}=0$. Since $0<r<n$, $p(s)=\det(\widetilde S-sI)$ is reducible in $\Q[s]$,
contradicting the fact that $M$ is simple. This concludes the proof of the lemma.
\end{proof}
\bigskip

\begin{lemma}\label{lemma:boundary2}
Two objects $M$ and $M'$ are isomorphic in $\C$ if and only if $\widehat M$ and $\widehat M'$ are isomorphic in $\D$.
\end{lemma}
\begin{proof}
Let $R\in GL_n(\Q)$ denote the matrix of an isomorphism between $M$ and $M'$ in $\C$. Then, one easily checks that the
matrix $R^{\oplus 4}\in GL_{4n}(\Q)$ defines an isomorphism between $\widehat M$ and $\widehat M'$ in $\D$.
Conversely, let us assume that we have an isomorphism between $\widehat M$ and $\widehat M'$ in $\D$. In particular,
it is an isomorphism between the $\Q$--vector spaces $M^{\oplus 4}$ and $(M')^{\oplus 4}$. Let $H\in GL_{4n}(\Q)$ be a matrix of this
isomorphism with respect to bases of $M^{\oplus 4}$ and $(M')^{\oplus 4}$ given by four copies of some fixed bases of $M$ and $M'$.
Since $h\pi_i=\pi_ih$ for $i=1,2$, $H$ is necessarily of the form $H=\begin{pmatrix}P&0\cr 0&Q\end{pmatrix}$ for some
$P,Q\in GL_{2n}(\Q)$. Furthermore, $H$ also satisfies the equation $H\widehat S=\widehat S'H$. Here, $\widehat S$ (resp.$\widehat S'$)
is the matrix obtained via equation \eqref{hatS} from the matrix $S$ (resp. $S'$) giving the action of $s\in P_1=\Z[s]$ on $M$ (resp. $M'$).
This implies easily that $P=Q=\begin{pmatrix}R&0\cr 0&R\end{pmatrix}$ for some $R\in GL_n(\Q)$ that satisfies $RS=S'R$.
This matrix $R$ defines an isomorphism in $\C$ between $M$ and $M'$.
\end{proof}

If $M$ is a $(-1)$--self-dual object in $\C$ via $\phi\co M\to M^*$, then $\widehat M$ is a $(-1)$--self-dual object in $\D$
via $\what\phi\co\widehat M\to(\widehat M)^*$. Hence, Lemmas \ref{lemma:boundary1} and \ref{lemma:boundary2} imply that the
Bing doubling map $(M,\phi)\mapsto(\widehat M,\what \phi)$ induces an inclusion
$\overline{\M}^s(\C,-1)\hookrightarrow\overline{\M}^s(\D,-1)$. By Hermitian devissage, it induces
homomorphisms $\varphi_M$ for all $M\in\overline{\M}^s(\C,-1)$ which fit in the following commutative diagram:
\[
\xymatrix{
W^{-1}(\C)\ar[r]\ar[d]^\cong&W^{-1}(\D)\ar[d]^\cong\\
\bigoplus_{M}W^{-1}(\C|_M)\ar[r]\ar@{->>}[d]^{p_M}&\bigoplus_{N}W^{-1}(\D|_N)\ar@{->>}[d]^{p_{\wwhat M}}\\
W^{-1}(\C|_M)\ar[r]^{\varphi_M}&W^{-1}(\D|_{\widehat M})
}
\] 
We are left with the proof that $\varphi_M$ is injective for all $M\in\overline{\M}^s(\C,-1)$. Consider the additive functor
$F\co\C\to\D$ given by $F(M)=\widehat M$ and $F(\alpha\co M\to M')=(\alpha^{\oplus 4}\co \widehat M\to \widehat M')$.

\begin{lemma}\label{lemma:boundary3}
If $M$ be a simple object in $\C$, then the ring homomorphism $\End_\C M\to\End_\D\widehat M$ induced by $F$ is an isomorphism.
\end{lemma}
\begin{proof}
A simple representation over $\Q$ of $P_1=\Z[s]$ is of the form $K=\Q[s]/(p)$, where $p\in\Q[s]$ is an irreducible polynomial.
Therefore, $\End_\C M=\End_{\Q[s]}K=K$ is a field. This implies that the homomorphism induced by $F$ is injective. To check
that it is onto, fix an element
$\beta\in\End_\D\widehat M$. Since $\beta\pi_i=\pi_i\beta$ for $i=1,2$, a matrix for $\beta$ is necessarily of the form
$\begin{pmatrix}P&0\cr 0&Q\end{pmatrix}$ for some $P,Q\in\Mat_{2n}(\Q)$. As in the proof of \fullref{lemma:boundary2},
the equation $\beta s=s\beta$ easily implies that $P=Q=\begin{pmatrix}R&0\cr 0&R\end{pmatrix}$ for some $R\in\Mat_n(\Q)$
such that $RS=SR$. This matrix $R$ defines an element $\alpha\in\End_\C M$ such that $F(\alpha)=\beta$.
\end{proof}

In particular, $F$ induces an isomorphism $F_*\co W^1(\End_\C M)\to W^1(\End_\D\widehat M)$.
Let $(M,b)$ be a $(-1)$--hermitian form in $\C$, and let $(\widehat M,\what b)$ be the corresponding $(-1)$--hermitian form in $\D$,
where a matrix for $\what b$ is obtained from a matrix for $b$ via \eqref{hatT}. Using \fullref{lemma:boundary3},
it is routine to check that for all $M\in\overline{\M}^s(\C,-1)$, the following diagram commutes:
\[
\xymatrix{W^{-1}(\C|_M)\ar[d]^\cong_{\Theta_{M,b}}\ar[r]^{\varphi_M}&W^{-1}(\D|_{\widehat M})\ar[d]^\cong_{\Theta_{\widehat M,\what b}}\\
W^1(\End_\C M)\ar[r]^{F_*}_\cong &W^1(\End_\D\widehat M)
}
\]
This concludes the proof of \fullref{thm:boundary}.\qed

\section{The Rasmussen invariant}\label{section:s}

In a beautiful paper \cite{Ras}, J Rasmussen used Khovanov homology to define a combinatorial invariant of smooth concordance of knots.
Its extension to links was carried out by A Beliakova and S Wehrli \cite{B-W}. We shall now show that this invariant,
unlike the ones listed in \fullref{section:Bing}, can tell that some Bing doubles are not smoothly slice.

Recall that the  Rasmussen invariant assigns to each oriented link $L$ in $S^3$ an integer $s(L)\in\Z$.
It satisfies the following properties (see \cite{B-W}).
\begin{romanlist}
\item{If $C$ is a smooth oriented cobordism in $S^3\times[0,1]$ between two oriented links $L$ and $L'$ such that every component of $C$ meets $L$, then
\[
|s(L)-s(L')|\le -\chi(C).
\]}
\item{If $L$ is represented by a positive diagram with $n$ crossings and $k$ Seifert circles, then $s(L)=n-k+1$.}
\item{If $L$ has $m$ components, then $s(L)\equiv m-1\pmod{2}$.}
\end{romanlist}

Note that Property (i) implies that $s$ is an invariant of smooth concordance, and Property (ii) implies that the $s$--invariant of the $m$--component
unlink is $1-m$. In particular, if $L$ is a smoothly slice 2--component link, then $s(L)=-1$.

\begin{lemma}\label{lemma:sBing}
The  Rasmussen invariant of a Bing double is equal to $\pm 1$.
\end{lemma}
\parpic[r]{$\begin{array}{c}
\includegraphics[height=2.5cm]{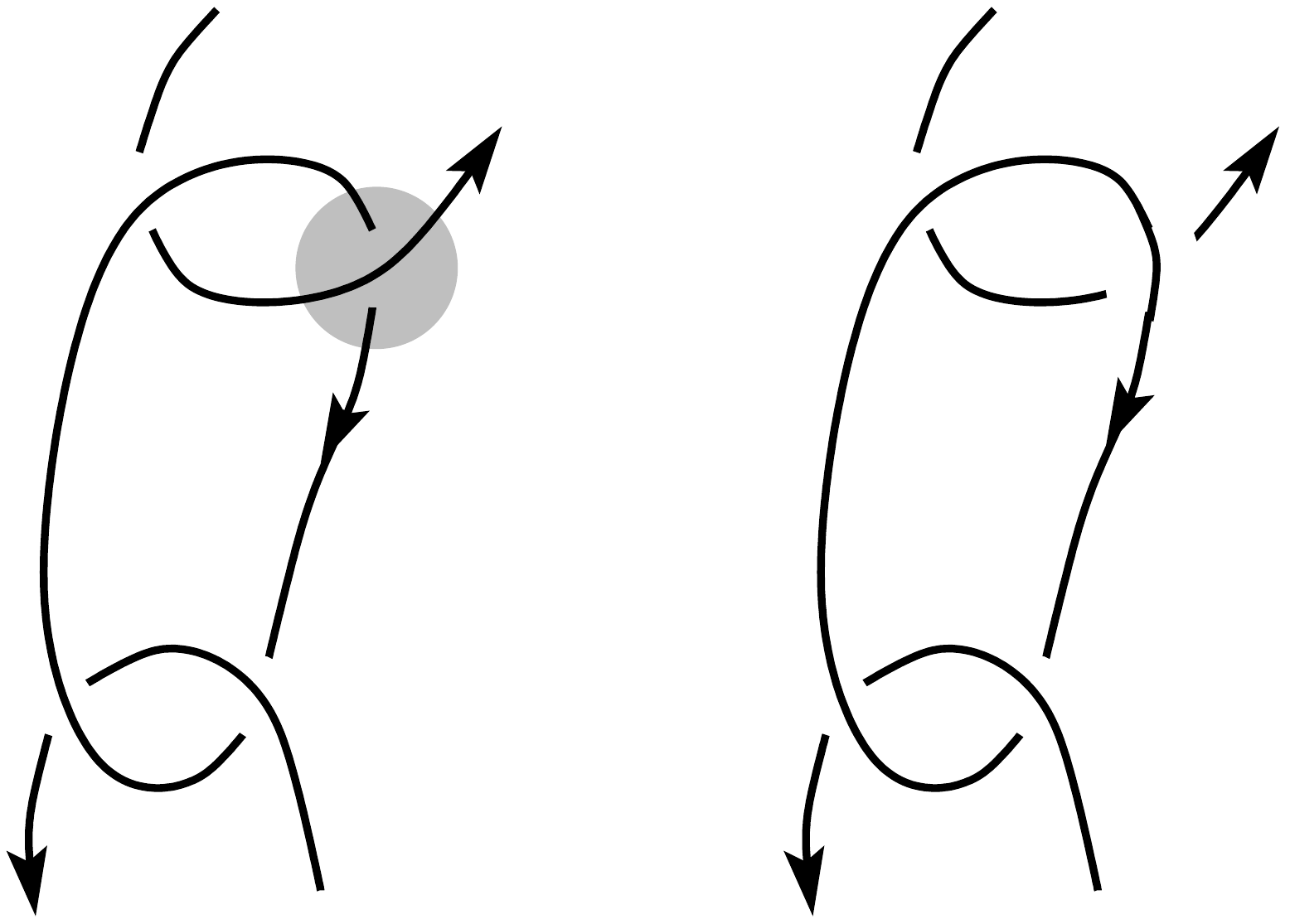}
\end{array}$}
\begin{proof}
The crossing change illustrated opposite turns any Bing double $B(K)$ into a positive Hopf link whose $s$--invariant is 1 by
Property (ii). This corresponds to a cobordism between $B(K)$ and the positive Hopf link with Euler characteristic $-2$.
By Property (i), $|s(B(K))-1|\le 2$. Similarly, the modification of one of the positive crossings in the diagram given here
turns $B(K)$ into a negative Hopf link, whose $s$--invariant is equal to $-1$. Hence, $|s(B(K))+1|\le 2$. The lemma now follows
from Property (iii).
\end{proof}

\begin{figure}[Htb]
\centerline{\psfig{file=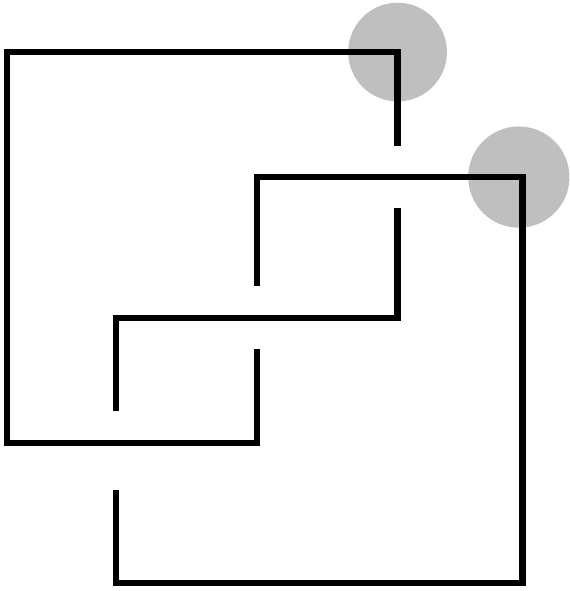,height=2.5cm}}
\caption{A diagram $D$ for the positive trefoil knot, with the 2 northeastern corners highlighted. On this example, $\mathrm{tb}(D)=1$.}
\label{fig:tb}
\end{figure}

Let us recall briefly the (combinatorial) definition of the Thurston--Bennequin invariant of a knot.
We refer to Rudolph \cite{Ru2} for the analytic definition and further details. Every knot admits a diagram $D$ consisting of vertical
and horizontal segments, the horizontal segment passing over the vertical one at each crossing.
Let $\mathrm{tb}(D)$ denote the writhe of $D$ minus the number of `northeastern corners' of $D$.
(See \fullref{fig:tb} for an example.) The {\em Thurston--Bennequin invariant\/} $\mathrm{TB}(K)$ of a knot $K$
is the maximum value of $\mathrm{tb}(D)$ over all such diagrams for $K$.

We shall be interested in knots $K$ such that $\mathrm{TB}(K)\ge 0$.
Let us mention that if a non-trivial knot $K$ is {\em strongly quasi-positive\/} (ie, if it is the closure of the product of braids of the form
$\sigma_{ij}=(\sigma_i\cdots\sigma_{j-2})\sigma_{j-1}(\sigma_i\cdots\sigma_{j-2})^{-1}$ with $1\le i<j<n$), then $\mathrm{TB}(K)\ge 0$.
In particular, all non-trivial positive knots (non-trivial knots that admit a diagram with only positive crossings) have non-negative Thurston--Bennequin invariant. 

\begin{prop}\label{prop:sTB}
If a knot $K$ has non-negative Thurston--Bennequin invariant, then $s(B(K))=1$.
\end{prop}

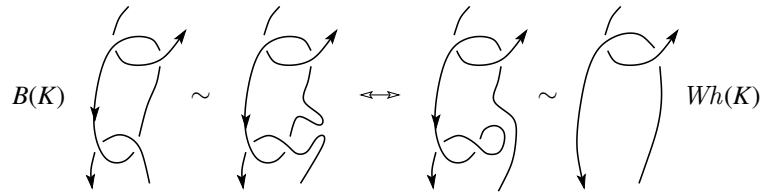
\begin{figure}[Htb]
\labellist\small\hair 2.5pt
\pinlabel {$B(K)$} at -80 165
\pinlabel {$\sim$} at 200 165
\pinlabel {$\sim$} at 800 165
\pinlabel {$\mathit{Wh}(K)$} at 1100 165
\endlabellist
\centerline{\psfig{file=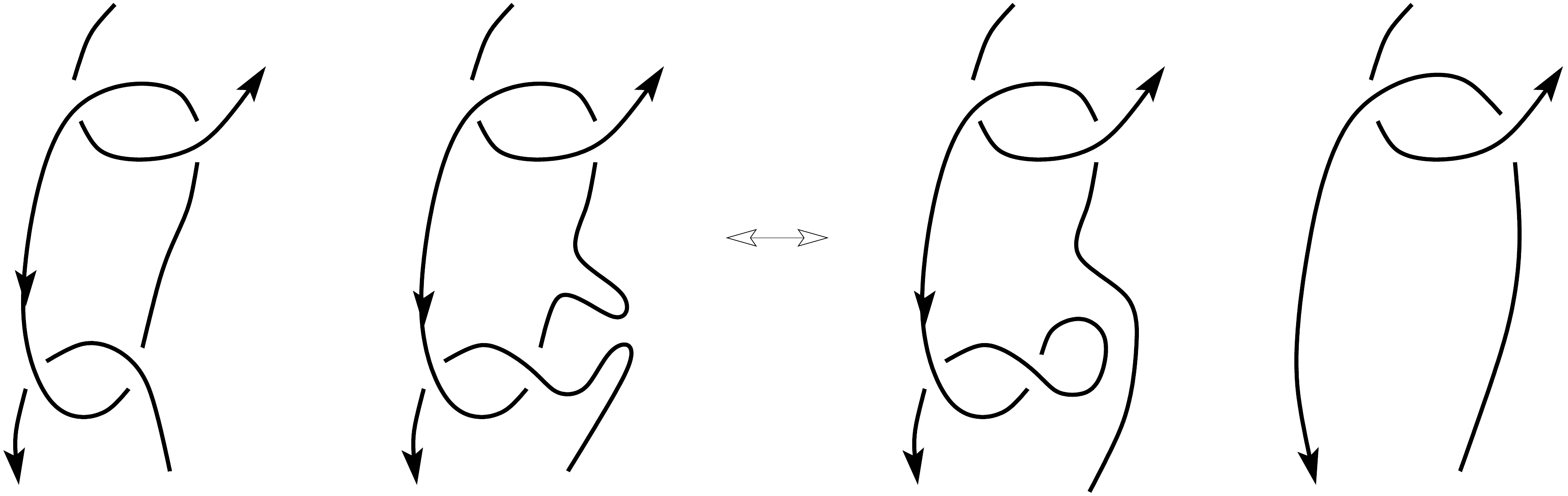,height=2.5cm}}
\caption{A cobordism between the Bing double and the untwisted positive Whitehead double of $K$}
\label{fig:BD}
\end{figure}

\begin{proof}
Consider the local transformation of the Bing double illustrated in \fullref{fig:BD}. It describes a cobordism $C$ between the Bing double $B(K)$
and the untwisted positive Whitehead double $\mathit{Wh}(K)$ of $K$ with $\chi(C)=-1$. Therefore, $|s(B(K))-s(\mathit{Wh}(K))|\le 1$.
Now, there is an obvious cobordism of Euler characteristic $-1$ between $\mathit{Wh}(K)$ and the positive Hopf link. Hence, $|s(\mathit{Wh}(K))-1|\le 1$.
Since $s(\mathit{Wh}(K))$ is even, it is equal to $0$ or $2$. Gathering these results with \fullref{lemma:sBing}, we conclude that the only possible values
of the pair $(s(B(K)),s(\mathit{Wh}(K)))$ are $(-1,0)$, $(1,0)$ and $(1,2)$. In particular, if $s(\mathit{Wh}(K))=2$, then $s(B(K))=1$.

Building on ideas of L Rudolph, C Livingston showed that if a knot satisfies $\mathrm{TB}(K)\ge 0$, then $s(\mathit{Wh}(K))=2$. This computation
was originally performed for the Ozsv\'ath--Szab\'o knot concordance invariant, but using only properties shared by the Rasmussen invariant
(see Livingston \cite[Theorem 12]{Liv} and Livingston--Naik \cite[Theorem 2]{L-N}). Therefore, if $\mathrm{TB}(K)\ge 0$, then $s(B(K))=1$.
\end{proof}

As an immediate consequence, we see that if $\mathrm{TB}(K)\ge 0$, then $B(K)$ is not smoothly slice. This statement is by no means new. Indeed,
L Rudolph \cite{Ru2} showed that if $\mathrm{TB}(K)\ge 0$, then $\mathit{Wh}(K)$ is not smoothly slice. By the cobordism illustrated in
\fullref{fig:BD}, this implies that $B(K)$ is not smoothly slice. However, it is interesting to see that this fact can be recovered using only
basic properties of the Rasmussen invariant for links.

\fullref{prop:sTB} also provides a wide class of links for which the Rasmussen invariant is a better obstruction to sliceness then the
multivariable Levine--Tristram signature (recall \fullref{section:Bing}). These examples are the first non-split links with
this property (compare \cite[Section 6.3]{B-W}). Since their components are trivial, our examples also show that the obstruction
to the sliceness of a link $L$ provided by $s(L)$ is more then just the obstruction to the sliceness of each components $L_i$ provided by $s(L_i)$. 

One easily checks that the existence of an invariant of links that satisfies Properties (i) and (ii) above implies the slice
Bennequin inequality. (The argument for knots given in Shumakovitch \cite{Shu} extends to links.) By \fullref{prop:sTB} and the end of \fullref{section:Bing}, we see
that the Rasmussen invariant is actually stronger than the slice Bennequin inequality for detecting links that are not smoothly slice.

Let us conclude this article with one last remark. Freedman showed that the untwisted Whitehead double of any knot $K$
is topologically slice (see Freedman and Quinn \cite{F-Q}).
By \fullref{cor:slice}, the link $B(\mathit{Wh}(K))$ is topologically slice. On the other hand, if $\mathrm{TB}(K)\ge 0$, then
$\mathrm{TB}(\mathit{Wh}(K))\ge 1$ by \cite[Proposition 3]{Ru2}. By \fullref{prop:sTB}, $B(\mathit{Wh}(K))$ is not smoothly slice.
Therefore, each knot $K$ with non-negative Thurston--Bennequin invariant induces a link $B(\mathit{Wh}(K))$ that is topologically but not smoothly slice, and whose components are trivial.

\medskip
{\bf Acknowledgements}\qua
The author wishes to express his thanks to Mathieu Baillif, Stefan Friedl, Lee Rudolph and Stefan Wehrli.
Above all, it is a pleasure to thank Peter Teichner for numerous discussions.
This work was supported by the Swiss National Science Foundation.

\bibliographystyle{gtart}
\bibliography{link}

\end{document}